\documentclass[12pt]{article}
\usepackage{mathrsfs}
\usepackage{epic,eepic,epsf,epsfig}
\usepackage{amsfonts,srcltx,mathrsfs}
\usepackage{multirow}
\usepackage{amsmath}
\usepackage{amssymb}
\usepackage{amsbsy}
\usepackage{graphicx}
\usepackage{amsfonts}
\usepackage{color}
\usepackage{setspace}

\newtheorem{lem}{Lemma}[section]%
\newtheorem{theorem}[lem]{Theorem}%

\def\nd{\mathrel{\bigm|\kern-.7em/}}

\def\f{\noindent}

\def\P\GammaL{\hbox{\rm P\GammaL}}

\begin{document}
\title{Eigenvalues, edge-disjoint perfect matchings and toughness of regular graphs}

\footnotetext{E-mails: zhangwq@pku.edu.cn}

\author{Wenqian Zhang\\
{\small School of Mathematics and Statistics, Shandong University of Technology}\\
{\small Zibo, Shandong 255000, P.R. China}}
\date{}
\maketitle

\begin{abstract}
Let $G$ be a connected $d$-regular graph of order $n$, where $d\geq3$. Let $\lambda_{2}(G)$ be the second largest eigenvalue of $G$. For even $n$, we show that $G$ contains $\left\lfloor\frac{2}{3}(d-\lambda_{2}(G))\right\rfloor$ edge-disjoint perfect matchings. This improves a result stated by Cioab\u{a},  Gregory and  Haemers \cite{CGH}. Let $t(G)$ be the toughness of $G$. When $G$ is non-bipartite, we give a sharp upper bound of $\lambda_{2}(G)$ to guarantee that $t(G)>1$. This enriches the previous results on this direction.

\bigskip

\f {\bf Keywords:} eigenvalue; edge-disjoint perfect matching; toughness; regular graph.

\end{abstract}

 \baselineskip 17 pt

\section{Introduction}

In this paper, all graphs considered are finite, undirected and without loops or multiple edges.  For a graph $G$, let $\overline{G}$ denote the complement of $G$. The vertex set and edge set of $G$ are denoted by $V(G)$ and $E(G)$, respectively. Let $e(G)=|E(G)|$.  For two disjoint subsets $X$ and $Y$ of $V(G)$, let $E_{G}(X,Y)$ be the set of edges with one end vertex in $X$ and the other end vertex in $Y$. Let $e_{G}(X,Y)=|E_{G}(X,Y)|$. For a subset $S$ of $V(G)$, let $G[S]$ be the subgraph induced by $S$, and let $G-S=G[V(G)-S]$. We also denote $G-\left\{u\right\}$ by $G-u$. For a subset $E$ of $E(G)$, let $G-E$ be the graph obtained from $G$ by deleting the edges in $E$.   For an integer $n\geq1$, let $K_{n}$ be the complete graph of order $n$. For two vertex-disjoint graphs $G_{1}$ and $G_{2}$, let $G_{1}\cup G_{2}$ be the disjoint union of them. For any terminology used but not defined here, one may refer to \cite{BH}.

For an $n\times n$  square matrix $B$ with real eigenvalues, let $\lambda_{1}(B)\geq \lambda_{2}(B)\geq\cdots\geq \lambda_{n}(B)$ denote its eigenvalues. Let $G$ be a graph with vertices $u_{1},u_{2},\ldots,u_{n}$. The {\em adjacency matrix} $A(G)$ of $G$ is an $n\times n$ square matrix $(a_{ij})$, where $a_{ij}=1$ if $u_{i}$ is adjacent to $u_{j}$, and $a_{ij}=0$ otherwise. The eigenvalues of $G$ are the eigenvalues of its adjacency matrix $A(G)$, which are denoted by $\lambda_{1}(G)\geq \lambda_{2}(G)\geq\cdots\geq \lambda_{n}(G)$. As is well known, $\lambda_{1}(G)=d$, if $G$ is a $d$-regular graph.

Let $G$ be a graph.  The {\em edge-connectivity} of $G$ is the minimum number of edges whose deletion induces a non-connected graph. For an integer $\ell$, the graph $G$ is called $\ell$-edge-connected if its edge-connectivity is at least $\ell$.  A {\em matching} of $G$ is a set of disjoint edges of $G$. A matching $M$  of $G$ is called perfect, if each vertex of $G$ is incident with some edge in $M$. Graph $G$ is called {\em factor-critical}, if $G-u$ has a perfect matching for any $u\in V(G)$.

Let $G$ be a connected $d$-regular graph of even order. The relationship between eigenvalues and perfect matchings of $G$ was studied by many researchers.
Brouwer and Haemers \cite{BH1} first gave sufficient conditions for the existence of a perfect matching in $G$ in terms of $\lambda_{3}(G)$.
The result in  \cite{BH1} was further improved in \cite{C,CG1}. In 2009, Cioab\u{a},  Gregory and  Haemers \cite{CGH} obtained the following sharp bound of $\lambda_{3}(G)$ to guarantee that $G$ has a perfect matching.

\begin{theorem}{\rm (Cioab\u{a},  Gregory and  Haemers \cite{CGH})}\label{matching CGH}
Let $G$ be a connected $d$-regular graph of even order, where $d\geq3$. Let $\theta$ be the largest root of $x^{3}-x^{2}-6x+2=0$. If
\begin{equation}
\lambda_{3}(G)<\left\{
\begin{array}{ll}
\theta=2.85577... & for~d=3,\\
\frac{1}{2}(d-2+\sqrt{d^{2}+12})& for~even~d\geq4,\\
\frac{1}{2}(d-3+\sqrt{(d+1)^{2}+16})& for~odd~d\geq5,
\end{array}
\right.\notag
\end{equation}
then $G$ has a perfect matching.
\end{theorem}

 As a corollary of Theorem \ref{matching CGH},  Cioab\u{a},  Gregory and  Haemers \cite{CGH} (see Corollary 4 therein) obtained the following result
 on the number of edge-disjoint perfect matchings in a regular graph of even order. This result was first stated in \cite{BH1}
in terms of Laplacian eigenvalues.

\begin{theorem} {\rm (Cioab\u{a},  Gregory and  Haemers \cite{CGH})}\label{disjoint matching}
Let $G$ be a connected $d$-regular graph of even order, where $d\geq3$. Then $G$ contains $\left\lfloor\frac{1}{2}(d-\lambda_{2}(G))+\frac{1}{2}\right\rfloor$ edge-disjoint perfect matchings.
\end{theorem}

Inspired by Theorem \ref{disjoint matching}, it is natural to ask the largest number $C$, such that a $d$-regular graph $G$ of even order contains $\left\lfloor C(d-\lambda_{2}(G))\right\rfloor$ edge-disjoint perfect matchings. Clearly, $C<1$. Theorem \ref{disjoint matching} states that $C\geq\frac{1}{2}$. The following Theorem \ref{main1} improves this to $C\geq\frac{2}{3}$.

\begin{theorem}\label{main1}
Let $G$ be a connected $d$-regular graph of even order, where $d\geq3$. Then $G$ contains $\left\lfloor\frac{2}{3}(d-\lambda_{2}(G))\right\rfloor$ edge-disjoint perfect matchings.
\end{theorem}

 For a graph $G$, let $c(G)$ be the number of components of $G$. The {\em toughness} $t(G)$ of $G$ is defined as $\min\left\{\frac{|S|}{c(G-S)}\right\}$, where the minimum is taken over all proper subsets $S\subset V(G)$ such that $c(G-S)\geq2$.
 Toughness is a measure of the connectivity of graphs. Let $G$ be a connected $d$-regular graph. The relationship between the toughness and eigenvalues of $G$ has been studied by many researchers. Let $\lambda=\max\left\{\lambda_{2}(G),-\lambda_{n}(G)\right\}$. Alon \cite{A} first showed that $t(G)>\frac{1}{3}(\frac{d^{2}}{d\lambda+\lambda^{2}}-1)$. Around the same time, Brouwer \cite{B} independently discovered a slightly better bound $t(G)>\frac{d}{\lambda}-2$. And in \cite{B1}, Brouwer conjectured that $t(G)\geq\frac{d}{\lambda}-1$. This conjecture was eventually settled by Gu \cite{G}. For the special case of toughness 1, Liu and Chen \cite{LC} improved the previous results as follows.

\begin{theorem}{\rm (Liu and Chen \cite{LC})}\label{toughness Liu}
Let $G$ be a connected $d$-regular graph with $d\geq3$. If
\begin{align*}
\begin{split}
\lambda_{2}(G)<\left\{
\begin{array}{lll}
d-1+\frac{3}{d+1} &for ~even~d,\\
\\
d-1+\frac{2}{d+1} &for ~odd~d,
\end{array}
\right.
\end{split}
\end{align*}
then $t(G)\geq1$.
\end{theorem}

Cioab\u{a} and Wong \cite{CW} further improved Theorem \ref{toughness Liu} by obtaining sharp upper bound of $\lambda_{2}(G)$ as in Theorem \ref{toughness CW}. Later, Cioab\u{a} and Gu \cite{CG2} established a spectral condition for $G$ with fixed edge-connectivity to
guarantee that $t(G)\geq1$. 

\begin{theorem}{\rm (Cioab\u{a} and Wong \cite{CW})}\label{toughness CW}
Let $G$ be a connected $d$-regular graph with $d\geq3$. If
\begin{align*}
\begin{split}
\lambda_{2}(G)<\left\{
\begin{array}{lll}
\frac{d-2+\sqrt{d^{2}+12}}{2} &for ~even~d,\\
\\
\frac{d-2+\sqrt{d^{2}+8}}{2} &for ~odd~d,
\end{array}
\right.
\end{split}
\end{align*}
then $t(G)\geq1$.
\end{theorem}

Let $G$ be a connected  $d$-regular graph with $d\geq3$. If $G$ is  bipartite, then $G$ is balanced. Clearly, we have $t(G)\leq1$ in this case. Thus, it is natural to ask the sharp upper bound of $\lambda_{2}(G)$  to guarantee that $t(G)>1$ when $G$ is non-bipartite. Our result on this direction is the following Theorem \ref{main2}.

For an integer $d\geq3$, let $G_{d}$ be the graph obtained from $\overline{K_{d}},\overline{K_{d-1}}$ and $K_{d}$ by connecting all the vertices in $V(\overline{K_{d}})$ to all the vertices in $V(\overline{K_{d-1}})$, and adding a perfect matching between $V(\overline{K_{d}})$ and $V(K_{d})$. Clearly, $G_{d}$ is a non-bipartite $d$-regular graph of order $3d-1$. Let $S=V(\overline{K_{d}})$. Then $G_{d}-S$ has exactly $|S|=d$ components, implying that $t(G_{d})\leq1$.

\begin{theorem}\label{main2}
Let $G$ be a connected non-bipartite $d$-regular graph with $d\geq3$. If $\lambda_{2}(G)\leq\frac{\sqrt{1+4(d-1)^{2}}-1}{2}$, then $t(G)>1$ unless $G=G_{d}$.
\end{theorem}

The rest of the paper is organized as follows. In Section 2, we include several known theorems to prove the main results of this paper. In Section 3, we give the proof Theorem \ref{main1}. The proof of Theorem \ref{main2} is given in Section 4.

\section{Preliminaries}

Let $G$ be a graph. An odd component of $G$ is a component of odd order. Let $o(G)$ denote the number of odd components of $G$. The following result was given by Tutte \cite{Tutte}.

\begin{theorem} {\rm (Tutte \cite{Tutte})}\label{Tutt}
A graph $G$ contains a perfect matching if and only if
$o(G-S)\leq|S|$ for any subset $S$ of $V(G)$.
\end{theorem}

Let $G$ be a graph of order $n$. For a partition $\left\{V_{1},V_{2},...,V_{m}\right\}$ of $V(G)$, the corresponding {\em quotient matrix} of the partition is an $m\times m$ matrix $(b_{ij})$, where $b_{ij}=\frac{e_{G}(V_{i},V_{j})}{|V_{i}|}$ for any $1\leq i\neq j\leq m$ and $b_{ii}=\frac{2e(G[V_{i}])}{|V_{i}|}$ for $1\leq i\leq m$. The partition is {\em equitable} if each vertex in $V_{i}$ has the same number of neighbors in $V_{j}$ for any $1\leq i, j\leq m$.

The following interlacing theorem can be found in \cite{BH} (Chapter 3), and is usually referred to as Cauchy eigenvalue interlacing.

\begin{theorem} {\rm (Brouwer and  Haemers \cite{BH})} \label{interlacing}
Let $G$ be a graph of order $n$, and let $B$ be the quotient matrix of a partition $\left\{V_{1},V_{2},...,V_{m}\right\}$ of $V(G)$. Then $\lambda_{i+n-m}(G)\leq \lambda_{i}(B)\leq\lambda_{i}(G)$ for any $1\leq i\leq m$. Moreover, if the partition is equitable, then all the eigenvalues (with multiplicities) of $B$ are also eigenvalues of $G$.
\end{theorem}

\section{Proof of Theorem \ref{main1}}

\begin{lem}\label{inequality}
Let $G$ be a graph and let $F$ be an induced subgraph of $G$. Assume that $V(F)$ can be partitioned in two parts $S$ and $T$, such that $\frac{2e(F[S])}{|S|}> a, \frac{2e(F[T])}{|T|}\geq b$ and  $\frac{e_{F}(S,T)}{|S|},\frac{e_{F}(S,T)}{|T|}\leq c$, where $a,b\geq0$ and $c>0$.  Then $\lambda_{2}(G)>\frac{a+b-\sqrt{(a-b)^{2}+4c^{2}}}{2}$.
\end{lem}

\f{\bf Proof:} By eigenvalue interlacing, we have $\lambda_{2}(G)\geq\lambda_{2}(F)$, since $F$ is an induced subgraph of $G$. Recall that $V(F)=S\cup T$. Let $d_{1}=\frac{2e(F[S])}{|S|},d_{2}=\frac{2e(F[T])}{|T|}$ and $d_{3}=\frac{e_{F}(S,T)}{|S|},d_{4}=\frac{e_{F}(S,T)}{|T|}$. Then $d_{1}>a,d_{2}\geq b$ and $d_{3},d_{4}\leq c$.
The quotient matrix $B$ of the partition $\left\{S, T\right\}$ is equal to
\begin{center}
$\left(\begin{array}{cc}
 d_{1}&d_{3}\\
d_{4}&d_{2}
\end{array}\right)$.
   \end{center}
By Theorem \ref{interlacing}, we have $\lambda_{2}(F)\geq\lambda_{2}(B)$.
By a simple calculation we have
 $$\lambda_{2}(B)=\frac{d_{1}+d_{2}-\sqrt{(d_{1}-d_{2})^{2}+4d_{3}d_{4}}}{2}\geq\frac{d_{1}+d_{2}-\sqrt{(d_{1}-d_{2})^{2}+4c^{2}}}{2}.$$

Let $f(x)=x+d_{2}-\sqrt{(x-d_{2})^{2}+4c^{2}}$, where $x\geq0$. By a simple calculation we have $$f'(x)=\frac{\sqrt{(x-d_{2})^{2}+4c^{2}}-(x-d_{2})}{\sqrt{(x-d_{2})^{2}+4c^{2}}}>0$$
 as $c>0$. Thus, the expression $\frac{d_{1}+d_{2}-\sqrt{(d_{1}-d_{2})^{2}+4c^{2}}}{2}$ is strictly increasing with respect to $d_{1}\geq0$. By symmetry, the expression $\frac{d_{1}+d_{2}-\sqrt{(d_{1}-d_{2})^{2}+4c^{2}}}{2}$ is also strictly increasing with respect to $d_{2}\geq0$. Recall that $d_{1}>a$ and $d_{2}\geq b$. Then $$\lambda_{2}(B)\geq\frac{d_{1}+d_{2}-\sqrt{(d_{1}-d_{2})^{2}+4c^{2}}}{2}>\frac{a+b-\sqrt{(a-b)^{2}+4c^{2}}}{2}.$$ It follows that $\lambda_{2}(G)\geq\lambda_{2}(F)\geq\lambda_{2}(B)>\frac{a+b-\sqrt{(a-b)^{2}+4c^{2}}}{2}$.
  This completes the proof. \hfill$\Box$

\medskip

\f{\bf Proof of Theorem \ref{main1}}. Let $n$ be the order of $G$, and let $\lambda_{2}=\lambda_{2}(G)$. If $d-\lambda_{2}<\frac{3}{2}$, then $\left\lfloor\frac{2}{3}(d-\lambda_{2})\right\rfloor=0$ and the theorem holds trivially. If $\frac{3}{2}\leq d-\lambda_{2}<3$, then $\lambda_{2}\leq d-\frac{3}{2}$ and $\left\lfloor\frac{2}{3}(d-\lambda_{2})\right\rfloor=1$. By Theorem \ref{matching CGH} (or Theorem \ref{disjoint matching}), we have that $G$ contains a perfect matching if $\lambda_{2}\leq d-1$. Thus, the theorem also holds in this case.
From now on, assume that $d-\lambda_{2}\geq3$. That is, $\lambda_{2}\leq d-3$. 

\medskip

\f{\bf Claim 1}. $G$ is $d$-edge-connected.

\medskip

\f{\bf Proof of Claim 1}. Suppose that $G$ is not $d$-edge-connected. Then there is a proper subset $U$ of $V(G)$, such that $e_{G}(U,V(G)-U)\leq d-1$.
If $|U|\leq d$, then by regularity, $e_{G}(U,V(G)-U)\geq (d+1-|U|)|U|\geq d$, a contradiction. Hence $|U|\geq d+1$ and $|V(G)-U|\geq d+1$ similarly. Now partition $V(G)$ in two parts $U$ and $V(G)-U$. Its quotient matrix of this partition is easy to calculate. By Theorem \ref{interlacing}, we have $\lambda_{2}\geq d-\frac{d-1}{d+1}-\frac{d-1}{d+1}>d-2$. But this contradicts that $\lambda_{2}\leq d-3$. Hence $G$ is $d$-edge-connected. This finishes the proof Claim 1. \hfill$\Box$

Let $t$ be the maximum integer such that $G$ contains $t$ edge-disjoint perfect matchings. By Theorem \ref{matching CGH}, we see that $t\geq1$ as $\lambda_{2}\leq d-3$. If $t\geq\left\lfloor\frac{2}{3}(d-\lambda_{2})\right\rfloor$, then the theorem holds. Thus we can assume that $t\leq\left\lfloor\frac{2}{3}(d-\lambda_{2})\right\rfloor-1$ in the following discussion. Let $M_{1},M_{2},...,M_{t}$ be $t$ edge-disjoint perfect matchings in $G$. Let $H$ be the spanning subgraph obtained from $G$ by deleting the edges in $\cup_{1\leq i\leq t}M_{i}$. Clearly, $H$ is a $(d-t)$-regular graph of order $n$. By the choice of $t$, $H$ has no perfect matchings. We will obtain a contradiction by the following two cases.

\medskip

\f{\bf Case 1}. $H$ is connected.

Since $H$ has no perfect matchings, by Theorem \ref{Tutt} there is a subset $S$ of $H$ such that $o(H-S)>|S|$. By parity we have $o(H-S)\geq|S|+2$ as $n$ is even. Let $q=o(H-S)$ and $s=|S|$. Then $q\geq s+2$ and $s\geq1$ as $H$ is connected. Let $Q_{1},Q_{2},...,Q_{q}$ be the odd components of $H-S$. Set $n_{i}=|Q_{i}|$ for any $1\leq i\leq q$. If there are at most $2$ odd components of $H-S$, say $Q_{1}$ and $Q_{2}$, such that $e_{H}(V(Q_{i}),V(G)-V(Q_{i}))<d-t$ for $1\leq i\leq2$, then by regularity we have
$$(d-t)s\geq\sum_{1\leq i\leq q}e_{H}(V(Q_{i}),V(G)-V(Q_{i}))\geq(d-t)(q-2)+2>(d-t)s,$$
a contradiction. Hence there are at least $3$ odd components of $H-S$, say $Q_{1},Q_{2}$ and $Q_{3}$, such that $e_{H}(V(Q_{i}),V(G)-V(Q_{i}))<d-t$ for $1\leq i\leq3$. Note that $H$ is $(d-t)$-regular. Similar to Claim 1, we can show that $n_{i}=|Q_{i}|\geq d-t+1$ for any $1\leq i\leq 3$. 
It follows that 
$$\frac{2e(G[V(Q_{i})])}{n_{i}}\geq \frac{2e(H[V(Q_{i})])}{n_{i}}=d-t-\frac{e_{H}(V(Q_{i}),V(G)-V(Q_{i}))}{n_{i}}>d-t-1$$
 for $1\leq i\leq3$.

Without loss of generality, assume that $n_{1}\leq n_{2}\leq n_{3}$. Recall that $H=G-\cup_{1\leq i\leq t}M_{i}$. Clearly, $$e_{G}(V(Q_{1}),V(Q_{2}))+e_{G}(V(Q_{1}),V(Q_{3}))\leq tn_{1}.$$
 Without loss of generality, assume that $e_{G}(V(Q_{1}),V(Q_{2}))\leq e_{G}(V(Q_{1}),V(Q_{3}))$. Then $e_{G}(V(Q_{1}),V(Q_{2}))\leq\frac{1}{2}tn_{1}$. 
 It follows that $$\frac{e_{G}(V(Q_{1}),V(Q_{2}))}{n_{1}}\leq\frac{1}{2}t$$
  and $$\frac{e_{G}(V(Q_{1}),V(Q_{2}))}{n_{2}}\leq\frac{1}{2}t.$$
Let $F=G[V(Q_{1})\cup V(Q_{2})]$. Then $\frac{2e(F[V(Q_{i})])}{n_{i}}>d-t-1$  and $\frac{e_{F}(V(Q_{1}),V(Q_{2}))}{n_{i}}\leq\frac{1}{2}t$ for $1\leq i\leq 2$. By Lemma \ref{inequality} (letting $a=b=d-t-1,c=\frac{1}{2}t$), we have $\lambda_{2}>d-t-1-\frac{1}{2}t$ and thus $t>\frac{2}{3}(d-\lambda_{2}-1)$. Since $t$ is an integer,  we obtain that $t\geq\left\lfloor\frac{2}{3}(d-\lambda_{2}-1)\right\rfloor+1=\left\lfloor\frac{2}{3}(d-\lambda_{2})+\frac{1}{3}\right\rfloor$. But this contradicts the assumption $t\leq\left\lfloor\frac{2}{3}(d-\lambda_{2})\right\rfloor-1$.

\medskip

\f{\bf Case 2}. $H$ is not connected.

Let $H_{1},H_{2},...,H_{m}$ be all the components of $H$, where $m\geq2$. Then $H_{i}$ is $(d-t)$-regular for each $1\leq i\leq m$. Let $h_{i}=|H_{i}|$ for $1\leq i\leq m$.

\f{\bf Subcase 2.1}. $m\geq3$.

Without loss of generality, assume that $h_{1}\leq h_{2}\leq h_{3}$. Recall that $H=G-\cup_{1\leq i\leq t}M_{i}$. Clearly, $$e_{G}(V(H_{1}),V(H_{2}))+e_{G}(V(H_{1}),V(H_{3}))\leq th_{1}.$$
 Without loss of generality, assume that $e_{G}(V(H_{1}),V(H_{2}))\leq e_{G}(V(H_{1}),V(H_{3}))$. Then $e_{G}(V(H_{1}),V(H_{2}))\leq\frac{1}{2}th_{1}$. It follows that $\frac{e_{G}(V(H_{1}),V(H_{2}))}{h_{i}}\leq\frac{1}{2}t$ for $1\leq i\leq2$.
Let $F=G[V(H_{1})\cup V(H_{2})]$. Then $\frac{2e(F[V(H_{i})])}{h_{i}}\geq d-t>d-t-1$ and $\frac{e_{F}(V(H_{1}),V(H_{2}))}{h_{i}}\leq\frac{1}{2}t$ for $1\leq i\leq 2$. By Lemma \ref{inequality} (letting $a=b=d-t-1,c=\frac{1}{2}t$), we have $\lambda_{2}>d-t-1-\frac{1}{2}t$. It follows that $t\geq\left\lfloor\frac{2}{3}(d-\lambda_{2}-1)\right\rfloor+1=\left\lfloor\frac{2}{3}(d-\lambda_{2})+\frac{1}{3}\right\rfloor$ as $t$ is an integer. But this contradicts the assumption $t\leq\left\lfloor\frac{2}{3}(d-\lambda_{2})\right\rfloor-1$.

\f{\bf Subcase 2.2}. $m=2$.

In this case, $H_{1}$ and $H_{2}$ are the only two components of $H$. Recall that $n_{i}=|H_{i}|$ for $1\leq i\leq2$. Since $n=h_{1}+h_{2}$, we see that $h_{1}$ and $h_{2}$ have the same parity.

\f{\bf Subcase 2.2.1}. $h_{1}$ and $h_{2}$  are even integers.

Since $H$ has no perfect matchings, without loss of generality, we can assume that $H_{1}$ has no perfect matchings.  By Theorem \ref{Tutt} there is a subset $S$ of $H_{1}$ such that $o(H_{1}-S)>|S|$. By parity we have $o(H_{1}-S)\geq|S|+2$ as $h_{1}$ is even. By a similar discussion to $H_{1}$ as $H$ in Case 1, we can obtain three odd components of $H_{1}-S$, say $Q_{1},Q_{2},Q_{3}$, such that they satisfy the same conditions as in Case 1. In a similar way, we can obtain a contradiction.

\f{\bf Subcase 2.2.2}. $h_{1}$ and $h_{2}$  are odd integers. In this case, there are still two subcases to consider.

\f{\bf Subcase 2.2.2.1}. At least one of $H_{1}$ and $H_{2}$ is not factor-critical.

 Without loss of generality, assume that $H_{1}$ is not factor-critical. Then there is a vertex $u$ of $H_{1}$, such that $H_{1}-u$ has no perfect matchings. By Theorem \ref{Tutt}, there is a subset $S_{0}$ of $V(H_{1}-u)$ such that $o(H_{1}-u-S_{0})>|S_{0}|$. Let $S=S_{0}\cup\left\{u\right\}$. Then $o(H_{1}-S)\geq|S|$. By parity we have $o(H_{1}-S)\geq|S|+1$ as $h_{1}$ is odd. Similar to Case 1, we can obtain two odd components of $H_{1}-S$, say $Q_{1}$ and $Q_{2}$, such that $e_{H_{1}}(V(Q_{i}),V(G)-V(Q_{i}))<d-t$ for $1\leq i\leq2$. Moreover, $|Q_{i}|\geq d-t+1$ and $\frac{2e(Q_{i})}{|Q_{i}|}>d-t-1$ for $1\leq i\leq2$. Let $Q_{3}=H_{2}$. Then  $\frac{2e(Q_{3})}{|Q_{3}|}=d-t>d-t-1$. Let $n_{i}=|Q_{i}|$ for $1\leq i\leq3$. Without loss of generality, assume that $n_{1}\leq n_{2}\leq n_{3}$. Recall that $H=H_{1}\cup H_{2}=G-\cup_{1\leq i\leq t}M_{i}$. Clearly, $$e_{G}(V(Q_{1}),V(Q_{2}))+e_{G}(V(Q_{1}),V(Q_{3}))\leq tn_{1}.$$
 Without loss of generality, assume that $e_{G}(V(Q_{1}),V(Q_{2}))\leq e_{G}(V(Q_{1}),V(Q_{3}))$. Then $e_{G}(V(Q_{1}),V(Q_{2}))\leq\frac{1}{2}tn_{1}$.
Let $F=G[V(Q_{1})\cup V(Q_{2})]$. Then $\frac{2e(F[V(Q_{i})])}{n_{i}}>d-t-1$  and $\frac{e_{F}(V(Q_{1}),V(Q_{2}))}{n_{i}}\leq\frac{1}{2}t$ for $1\leq i\leq 2$. By Lemma \ref{inequality} (letting $a=b=d-t-1,c=\frac{1}{2}t$), we have $\lambda_{2}>d-t-1-\frac{1}{2}t$. It follows that  $t\geq\left\lfloor\frac{2}{3}(d-\lambda_{2}-1)\right\rfloor+1=\left\lfloor\frac{2}{3}(d-\lambda_{2})+\frac{1}{3}\right\rfloor$. But this contradicts the assumption $t\leq\left\lfloor\frac{2}{3}(d-\lambda_{2})\right\rfloor-1$.

\f{\bf Subcase 2.2.2.2}. Both $H_{1}$ and $H_{2}$ are factor-critical.

Denote $X=V(H_{1})$ and $Y=V(H_{2})$. Then $|X|+|Y|=n$. Since $G$ is $d$-edge-connected by Claim 1, we have $e_{G}(X,Y)\geq d>t$. It follows that there is at least one of $M_{1},M_{2},...,M_{t}$, say $M_{t}$, such that $|M_{t}\cap E_{G}(X,Y)|\geq2$. Suppose that $uv$ and $xy$ are two edges in $M_{t}$ such that $u,x\in X$ and $v,y\in Y$. Since both $H[X]$ and $H[Y]$ are factor-critical, we have that $H[X]-u$ has a perfect matching $N_{1}$ and $H[Y]-v$ has a perfect matching $N_{2}$. Clearly, the union $N$ of $\left\{uv\right\},N_{1}$ and $N_{2}$ is a perfect matching of $G-\cup_{1\leq i\leq t-1}M_{i}$. Now $N,M_{1},M_{2},...,M_{t-1}$ are $t$ edge-disjoint perfect matchings of $G$. Let $K=G-N\cup(\cup_{1\leq i\leq t-1}M_{i})$. Note that $xy$ is an edge of $K$, and $uv$ is a non-edge of $K$. By the choice of $t$, we see that $K$ has no perfect matchings.

If both $K[X]$ and $K[Y]$ are connected, then $K$ is connected as $xy$ is an edge of $K$. We can obtain a contradiction in a similar way as Case 1. Now, without loss of generality, we can assume that $K[X]$ is not connected. Let $Q_{1}$  be a component of $K[X]$, and let $Q_{2}$ be the union of other components of $K[X]$. Let $Q_{3}=K[Y]$.
Note that $K$ is $(d-t)$-regular. Moreover, $e_{K}(V(Q_{1}),V(Q_{2}))=0$.
 Without loss of generality, assume that $u\in V(Q_{1})$. Recall that $uv$ is a non-edge of $K$ and $v\in V(Q_{3})$. Then $e_{K}(V(Q_{1}),Y)<|Q_{1}|,e_{K}(V(Q_{2}),Y)\leq|Q_{2}|$ and $|e_{K}(V(Q_{3}),X)<|Q_{3}|$. It follows that $\frac{2e(Q_{1})}{n_{1}},\frac{2e(Q_{3})}{n_{3}}>d-t-1$ and $\frac{2e(Q_{2})}{n_{2}}\geq d-t-1$.  Let $n_{i}=|Q_{i}|$ for $1\leq i\leq3$, and
 assume that $n_{1}\leq n_{2}\leq n_{3}$ without loss of generality. Recall that $K=G-N\cup(\cup_{1\leq i\leq t-1}M_{i})$. Clearly, $$e_{G}(V(Q_{1}),V(Q_{2}))+e_{G}(V(Q_{1}),V(Q_{3}))\leq (t+1)n_{1}.$$
 Without loss of generality, assume that $e_{G}(V(Q_{1}),V(Q_{2}))\leq e_{G}(V(Q_{1}),V(Q_{3}))$. Then $e_{G}(V(Q_{1}),V(Q_{2}))\leq\frac{1}{2}(t+1)n_{1}$.
Let $F=G[V(Q_{1})\cup V(Q_{2})]$. Then 
$$\frac{2e(F[V(Q_{1})])}{n_{1}}>d-t-1,\frac{2e(F[V(Q_{2})])}{n_{2}}\geq d-t-1$$
  and $\frac{e_{F}(V(Q_{1}),V(Q_{2}))}{n_{i}}\leq\frac{1}{2}(t+1)$ for $1\leq i\leq 2$. By Lemma \ref{inequality} (letting $a=b=d-t-1,c=\frac{1}{2}(t+1)$), we have $\lambda_{2}>d-t-1-\frac{1}{2}(t+1)$. It follows that  $t>\frac{2}{3}(d-\lambda_{2})-1$. This implies that $t\geq\left\lfloor\frac{2}{3}(d-\lambda_{2})-1\right\rfloor+1=\left\lfloor\frac{2}{3}(d-\lambda_{2})\right\rfloor$ as $t$ is an integer.  But this contradicts the assumption $t\leq\left\lfloor\frac{2}{3}(d-\lambda_{2})\right\rfloor-1$. This completes the proof. \hfill$\Box$

\section{Proof of Theorem \ref{main2}}

 Recall that $G_{d}$ is the graph obtained from $\overline{K_{d}},\overline{K_{d-1}}$ and $K_{d}$ by connecting all the vertices in $V(\overline{K_{d}})$ to all the vertices in $V(\overline{K_{d-1}})$, and adding a perfect matching between $V(\overline{K_{d}})$ and $V(K_{d})$. As is shown in the introduction, $G_{d}$ is a non-bipartite $d$-regular graph with $t(G_{d})\leq1$.
Clearly, $G_{d}$ has an equitable partition with three parts $V(\overline{K_{d}}),V(\overline{K_{d-1}})$ and $V(K_{d})$ . Very similar to the proof of $(i)$ of Theorem 3.1 in \cite{Z}, we can obtain the following Lemma \ref{tough lem}.

\begin{lem}\label{tough lem}
For $d\geq3$, let $G_{d}$ be the graph defined as above. Then $\lambda_{2}(G_{d})=\frac{\sqrt{1+4(d-1)^{2}}-1}{2}$. Moreover, $d-\frac{3}{2}+\frac{1}{8d}<\frac{\sqrt{1+4(d-1)^{2}}-1}{2}<d-\frac{3}{2}+\frac{3}{16d}$.
 \end{lem}

\f{\bf Proof of Theorem \ref{main2}}. We prove the theorem by contradiction. Suppose that $G\neq G_{d}$ is satisfying that $\lambda_{2}(G)\leq\frac{\sqrt{1+4(d-1)^{2}}-1}{2}$ and $t(G)\leq1$. We will obtain a contradiction. By the definition of toughness, there is a subset $S$ of $G$ such that $G-S$ is disconnected and $c(G-S)\geq |S|$. Let $q=c(G-S),s=|S|$, and let $Q_{1},Q_{2},...,Q_{q}$ be the components of $G-S$. Then $q\geq s$ and $q\geq2$.

\f{\bf Claim 1}. $q=s$.

\f{\bf Proof of Claim 1}. If $q>s$, then by a similar discussion as in Theorem \ref{main1}, we can show that there are two components of $G-S$, say $Q_{1}$ and $Q_{2}$, such that $e_{G}(V(Q_{i}),V(G)-V(Q_{i}))<d$ for $1\leq i\leq2$. It follows that $|Q_{i}|>d$ and $\frac{2e(Q_{i})}{|Q_{i}|}>d-1$ for $1\leq i\leq2$. Note that $e_{G}(V(Q_{1}),V(Q_{2}))=0$. By eigenvalue interlacing, we have $$\lambda_{2}(G)\geq\min\left\{\lambda_{1}(Q_{1}),\lambda_{1}(Q_{2})\right\}\geq
\min\left\{\frac{2e(Q_{1})}{|Q_{1}|},\frac{2e(Q_{2})}{|Q_{2}|}\right\}>d-1.$$
 This contradicts that $\lambda_{2}(G)\leq\frac{\sqrt{1+4(d-1)^{2}}-1}{2}<d-1$. This finishes the proof of Claim 1. \hfill$\Box$

Since $q=s$ and $G$ is not bipartite, by regularity, there is at least one non-singleton component of $G-S$. Thus there are only two cases to consider as follows.

\f{\bf Case 1}. There are at least two components of $G-S$ which are not singleton components.

We can find two non-singleton components of $G-S$, say $P$ and $Q$, such that $e_{G}(V(P),S)+e_{G}(V(Q),S)\leq 2d$ by regularity. Without loss of generality, assume that $e_{G}(V(P),S)\leq e_{G}(V(Q),S)$. Then $e_{G}(V(P),S)\leq d$, implying that $|P|\geq d$.

Now we show that $|Q|\geq d$. If $|Q|<d$, then $2\leq|Q|\leq d-1$. Thus $e_{G}(V(Q),S)\geq |Q|(d+1-|Q|)\geq2d-2$. Since $e_{G}(V(P),S)+e_{G}(V(Q),S)\leq 2d$, we have $e_{G}(V(P),S)\leq2$. It follows that $|P|\geq d+1$ and $|V(G)|-|P|\geq d+1$. Partition $V(G)$ in two parts $V(P)$ and $V(G)-V(P)$. By Theorem \ref{interlacing}, we have $$\lambda_{2}(G)\geq d-\frac{4}{d+1}>d-\frac{3}{2}+\frac{3}{16d}>\frac{\sqrt{1+4(d-1)^{2}}-1}{2},$$
 a contradiction. Therefore, $|Q|\geq d$.

 Let $e_{G}(V(P),S)=x\leq d$. By Theorem \ref{interlacing}, we have
 $$\lambda_{2}(G)\geq \min\left\{\lambda_{1}(P),\lambda_{1}(Q)\right\}\geq \min\left\{d-\frac{x}{d},d-\frac{2d-x}{d}\right\}= d-\frac{2d-x}{d}.$$
  Partition $V(G)$ in two parts $V(P)$ and $V(G)-V(P)$. By Theorem \ref{interlacing}, we have $\lambda_{2}(G)\geq d-\frac{x}{d}-\frac{x}{d+2}$ as $|V(G)-V(P)|\geq d+2$. Therefore, we have
 $$\lambda_{2}(G)\geq \max\left\{d-\frac{2d-x}{d}, d-\frac{x}{d}-\frac{x}{d+2}\right\}\geq d-\frac{4d+4}{3d+4}$$
  (when $x=\frac{2d^{2}+4d}{3d+4}$). It is easy to check that $d-\frac{4d+4}{3d+4}\geq d-\frac{3}{2}+\frac{3}{16d}>\frac{\sqrt{1+4(d-1)^{2}}-1}{2}$, a contradiction.

\f{\bf Case 2}. There is exactly one non-singleton component of $G-S$, say $P$.

Note that $q\geq2$. There is at least one singleton component of $G-S$. It follows that $q=s\geq d$. By regularity we have $e_{G}(V(P),S)\leq d$, implying that $|P|\geq d$. Thus $n\geq 3d-1$. Set $e_{G}(V(P),S)=x$.

\f{\bf Subcase 2.1}.  $x<d$.

We obtain $x\leq d-2$ as $x=d-2e(G[S])$.
Partition $V(G)$ in two parts $V(P)$ and $V(G)-V(P)$. By Theorem \ref{interlacing}, we have
$$\lambda_{2}(G)\geq d-\frac{x}{|P|}-\frac{x}{n-|P|}\geq d-\frac{d-2}{d}-\frac{d-2}{2d-1},$$
 since $|P|,n-|P|\geq d$ and $n\geq 3d-1$. It is easy to check that $d-\frac{d-2}{d}-\frac{d-2}{2d-1}>d-\frac{3}{2}+\frac{3}{16d}>\frac{\sqrt{1+4(d-1)^{2}}-1}{2}$, a contradiction.

\f{\bf Subcase 2.2}. $x=d$.

If $n=3d-1$, then $|P|=d$ and $s=d$. By regularity, the graph $G$ is isomorphic to $G_{k}$, a contradiction.

It remains the case of $n\geq 3d$. Let $p=|P|\geq d$. Partition $V(G)$ in three parts $V(P),S$ and $V(G)-S\cup V(P)$.
The quotient matrix $B$ of this partition is equal to
\begin{center}
$\left(\begin{array}{ccc}
 d-\frac{d}{p}&\frac{d}{p}&0\\
\frac{d}{s}&0&d-\frac{d}{s}\\
0&d&0
\end{array}\right)$.
   \end{center}
By a calculation, we have $\lambda_{2}(B)=\frac{-\frac{d}{p}+\sqrt{(\frac{d}{p})^{2}+4d^{2}(1-\frac{1}{s})(1-\frac{1}{p})}}{2}$.
By Theorem \ref{interlacing}, we have
$$\lambda_{2}(G)\geq\lambda_{2}(B)=\frac{-\frac{d}{p}+\sqrt{(\frac{d}{p})^{2}+4d^{2}(1-\frac{1}{s})(1-\frac{1}{p})}}{2},$$
 where $p,s\geq d$ and $n\geq 3d$. By differentiation, it is easy to check that the expression $\frac{-\frac{d}{p}+\sqrt{(\frac{d}{p})^{2}+4d^{2}(1-\frac{1}{s})(1-\frac{1}{p})}}{2}$ is strictly increasing with respect to $p\geq d$ and also to $s\geq d$. Since $n\geq3d$, at least one of $p$ and $s$ is strictly larger than $d$. It follows that $$\frac{-\frac{d}{p}+\sqrt{(\frac{d}{p})^{2}+4d^{2}(1-\frac{1}{s})(1-\frac{1}{p})}}{2}>
 \frac{-\frac{d}{d}+\sqrt{(\frac{d}{d})^{2}+4d^{2}(1-\frac{1}{d})(1-\frac{1}{d})}}{2}=\frac{\sqrt{1+4(d-1)^{2}}-1}{2}.$$
 Thus $\lambda_{2}(G)>\frac{\sqrt{1+4(d-1)^{2}}-1}{2}$, a contradiction.
 This completes the proof. \hfill$\Box$

\medskip

\f{\bf Data availability statement}

\medskip

There is no associated data.

\medskip

\f{\bf Declaration of Interest Statement}

\medskip

There is no conflict of interest.

\end{document}